\documentclass[12pt,reqno]{amsart}
\usepackage{amsmath}
\usepackage{amsfonts}
\usepackage{amssymb}
\usepackage{amsthm}
\usepackage{amscd}
\usepackage{amsgen}
\usepackage{latexsym}
\usepackage{url}

\catcode`\@=11
\@namedef{subjclassname@2010}{%
\textup{2010} Mathematics Subject Classification}
\catcode`\@=12

\newtheorem{defn}{Definition}[section]
\newtheorem{thm}{Theorem}[section]
\newtheorem{prop}{Proposition}[section]

\newtheorem{lemm}{Lemma}[section]
\newtheorem{corr}{Corollary}[section]
\newtheorem{ex}{Example}[section]

\numberwithin{equation}{section}

\allowdisplaybreaks

\newcommand{\Ric}{\operatorname{Ric}}
\newcommand{\scal}{\operatorname{scal}}
\newcommand{\tr}{\operatorname{tr}}

\def\T{{\mathcal T}}            
\def\V{{\mathcal V}}            

\def\W{{\sf W}}                 
\def\B{{\mathcal B}}            
\def\Rho{{\sf P}}               
\def\QR{{\mathcal Q}}      

\def\J{{\sf J}}                 

\def\r{{\mathbb R}}             
\def\f{{\frac{n}{2}}}

\title{On Branson's $Q$-curvature of order eight}

\author{Andreas Juhl}

\address{Humboldt-Universit\"at, Institut f\"ur Mathematik, Unter den Linden,
D-10099 Berlin} \email{ajuhl@math.hu-berlin.de}

\begin{document}

\begin{abstract} We prove a universal recursive formulas for Branson's
$Q$-curvature of order eight in terms of lower-order $Q$-curvatures,
lower-order GJMS-operators and holographic coefficients. The results
prove a special case of a conjecture in \cite{juhl-power}.
\end{abstract}


\maketitle

\centerline \today

\tableofcontents

\footnotetext{The work was supported by SFB 647
``Space-Time-Matter'' of DFG.}

\section{Introduction and statement of results}\label{intro}

It is well-known that on any Riemannian manifold $(M,g)$ of
dimension $n \ge 2$, the second-order differential operator
\begin{equation}\label{yamabe}
P_2(g) = \Delta_g - \left(\f-1\right) \frac{\scal(g)}{2(n-1)}
\end{equation}
is conformally covariant in the sense that
$$
e^{(\f+1)\varphi} P_2(e^{2\varphi}g)(u) = P_2(g)(e^{(\f-1)\varphi}u)
$$
for all $\varphi \in C^\infty(M)$ and all $u\in C^\infty(M)$. Here,
$\Delta_g$ denotes the Laplace-Beltrami operator of the metric $g$
and $\scal(g)$ is the scalar curvature of $g$. The operator $P_2$ is
called the {\em conformal Laplacian} or {\em Yamabe operator}. More
generally, in \cite{GJMS} Graham et al. proved that on any
Riemannian manifold $(M,g)$ of {\em even} dimension $n$, there
exists a finite sequence $P_2(g), P_4(g), \dots, P_n(g)$ of
geometric differential operators of the form
$$
\Delta_g^N + \mbox{lower order terms}
$$
so that
$$
e^{(\f+N)\varphi} P_{2N}(e^{2\varphi}g)(u) =
P_{2N}(g)(e^{(\f-N)\varphi}u).
$$
The operators $P_{2N}(g)$ are {\em geometric} in the sense that the
lower order terms are completely determined by the metric and its
curvature. On the flat space $\r^n$, there is no non-trivial
curvature, and we have $P_{2N} = \Delta^N$. We shall follow common
practice by referring to these operators as to the GJMS-operators.

The constant terms of the GJMS-operators lead to the notion of
Branson's $Q$-curvature (see \cite{bran-2}). The critical
GJMS-operator $P_n$ is special in the sense that it has vanishing
constant term, that is $P_n(g)(1) = 0$. More generally, for $2N <
n$, it is natural to write the constant term of $P_{2N}$ in the form
\begin{equation}\label{Q-curv}
P_{2N}(g)(1) = (-1)^N \left(\f-N\right) Q_{2N}(g)
\end{equation}
with a scalar Riemannian curvature invariant $Q_{2N}(g) \in
C^\infty(M)$ of order $2N$. With this convention, the critical
$Q$-curvature $Q_n(g)$ can be defined through $Q_{2N}$, $2N<n$, by
continuation.\footnote{The signs in \eqref{Q-curv} are required by
the convention that $-\Delta$ is non-negative.}

Since the algorithmic definition in \cite{GJMS} is quite involved, a
direct derivation of formulas for $P_{2N}$ in terms of the metric is
very complicated if possible at all. However, in the simplest cases
$N=1$ and $N=2$ such evaluations are well-known to yield the
familiar Yamabe operator \eqref{yamabe} and the Paneitz-operator
\begin{equation}\label{pan-op}
P_4 = \Delta^2 + \delta ( (n-2) \J - 4\Rho ) d + \left(\f-2\right)
\left(\f \J^2 - 2|\Rho|^2 - \Delta \J \right).
\end{equation}
Here we use the notation
$$
\J = \frac{\scal}{2(n-1)} \quad \mbox{and} \quad \Rho =
\frac{1}{n-2}(\Ric - \J g).
$$
$\Rho$ is called the Schouten tensor. In \eqref{pan-op}, it is
regarded as an endomorphism of $\Omega^1(M)$. Eq.~\eqref{pan-op}
shows that, on manifolds of dimension $n \ge 4$,
\begin{equation}\label{q4-gen}
Q_4 = \f \J^2 - 2|\Rho|^2 - \Delta \J.
\end{equation}
In particular, on manifolds of dimension four, the critical
$Q$-curvature is given by
\begin{equation}\label{q4-crit}
Q_4 = 2(\J^2-|\Rho|^2) - \Delta \J.
\end{equation}
It has often been said that the complexity of GJMS-operators and
$Q$-curvatures increases exponentially with their order. It is
tempting to compare this with the complexity of heat coefficients.
Some explicit formulas for $Q_6$ and $Q_8$ in terms of the Schouten
tensor $\Rho$, the Weyl tensor $\W$ and their covariant derivatives
were derived in \cite{G-P}. The enormous complexity of these
formulas indicates that it is extremely hard to unveil the structure
of high order $Q$-curvatures. A crucial part of the problem is to
decide about the most natural way of stating the results.

In \cite{juhl-book}, we introduced and developed the idea to
investigate $Q$-curvature from a conformal submanifold perspective.
In particular, we introduced the notion of residue families
$D_{2N}^{res}(g;\lambda)$. These are certain families of local
operators which contain basic information on the structure of
$Q$-curvatures and GJMS-operators. Besides motivations by
representation theory, the approach builds on the interpretation of
GJMS-operators as residues of the scattering operator of conformally
compact Einstein metrics (see \cite{GZ}). The residue families
satisfy systems of recursive relations which can be used to reveal
the recursive structure of $Q$-curvatures and GJMS-operators. Along
such lines, we found recursive formulas for the critical
$Q$-curvatures $Q_6$ and $Q_8$, which express these quantities in
terms of respective lower order GJMS-operators and lower order
$Q$-curvatures. In in \cite{juhl-power} and \cite{FJ}, these methods
were further developed and led to the formulation of a number of
conjectures.

In \cite{juhl-book}, the discussion of recursive formulas for $Q_8$
was depending on some technical assumptions. Here we remove these
assumptions.

The formulation of the main results requires to define one more
ingredient. For a given metric $g$ on the manifold $M$ of even
dimension $n$, let
\begin{equation}\label{PE1}
g_+ = r^{-2} (dr^2 + g_r)
\end{equation}
with
\begin{equation}\label{PE2}
g_r = g + r^2 g_{(2)} + \dots + r^{n-2} g_{(n-2)} + r^n (g_{(n)}+
\log r \bar{g}_{(n)}) + \cdots
\end{equation}
be a metric on $M \times (0,\varepsilon)$ so that the tensor
$\Ric(g_+) + n g_+$ satisfies the Einstein condition
\begin{equation}\label{einstein}
\Ric(g_+) + n g_+ = O(r^{n-2})
\end{equation}
together with a certain vanishing trace condition. These conditions
uniquely determine the coefficients $g_{(2)}, \dots, g_{(n-2)}$.
They are given as polynomial formulas in terms of $g$, its inverse,
the curvature tensor of $g$, and its covariant derivatives. The
coefficient $\bar{g}_{(n)}$ and the quantity $\tr g_{(n)}$ are
determined as well. Moreover, $\bar{g}_{(n)}$ is trace-free, and the
trace-free part of $g_{(n)}$ is undetermined. A metric $g_+$ with
these properties is called a Poincar\'e-Einstein metric with
conformal infinity $[g]$. For full details see \cite{FG-final}.

The volume form of $g_+$ can be written as
$$
vol(g_+) = r^{-n-1} v(r) dr vol(g),
$$
where
$$
v(r) = vol(g_r)/vol(g) \in C^\infty(M).
$$
The coefficients $v_0,\dots,v_n$ in the Taylor series
$$
v(r) = v_0 + v_2 r^2 + v_4 r^4 + \cdots + v_n r^n + \cdots
$$
are known as the {\em renormalized volume} coefficients
(\cite{G-vol}, \cite{G-ext}) or {\em holographic} coefficients
(\cite{juhl-book}, \cite{BJ}). The coefficient $v_{2j} \in
C^\infty(M)$ is given by a local formula which involves at most $2j$
derivatives of the metric. Note also that $v_n$ is uniquely
determined by $g$ since $\tr g_{(n)}$ is uniquely determined by $g$.
It is called the holographic anomaly. Explicit formulas for the
holographic coefficients $v_2, v_4, v_6, v_8$ were derived in
\cite{G-ext}.

The first main result describes the {\em critical} $Q$-curvature of
order eight.

\begin{thm}\label{main} On manifolds of dimension $8$, Branson's
$Q$-curvature $Q_8$ is given by the formula
\begin{multline}\label{Q8-main}
Q_8 = -3 P_2(Q_6) - 3 P_6(Q_2) + 9 P_4(Q_4) \\
+ 8 P_2 P_4(Q_2) - 12 P_2^2 (Q_4) + 12 P_4 P_2 (Q_2) - 18 P_2^3(Q_2)
+ 3! 4! 2^8 w_8,
\end{multline}
where $w_8$ is the coefficient of $r^8$ in the Taylor series of
$\sqrt{v(r)}$.
\end{thm}

In terms of holographic coefficients, the quantity $w_8$ can be
expressed as
\begin{equation}\label{v8w8}
128 w_8 = 64 v_8 - 32 v_6 v_2 - 16 v_4^2 + 24 v_2^2 v_4 - 5 v_2^4
\end{equation}
(see Lemma \ref{sroot}).

A version of Theorem \ref{main} was proved in Section 6.13 of
\cite{juhl-book} under the assumption that the polynomial
$\V_8(\lambda)$ (see \eqref{V8}) vanishes. In the present paper, we
show that this assumption is vacuous (Proposition \ref{V8-van}).

Theorem \ref{main} confirms the special case $n=8$ and $N=4$ of a
conjectural formula for all $Q$-curvatures $Q_{2N}$ formulated in
\cite{juhl-power}. In connection with this conjecture, it is
important to recognize that the coefficients in \eqref{Q8-main} have
a uniform definition. In order to describe this, we introduce some
notation. A sequence $I=(I_1,\dots,I_r)$ of integers $I_j \ge 1$
will be regarded as a composition of the sum $|I| = I_1 + I_2 +
\cdots + I_r$, where two representations which contain the same
summands but differ in the order of the summands are regarded as
different. $|I|$ is called the size of $I$. For $I=(I_1,\dots,I_r)$,
we set
$$
P_{2I} = P_{2I_1} \circ \cdots \circ P_{2I_r}.
$$
For any composition $I$, we define the multiplicity $m_I$ by
\begin{equation}\label{m-form}
m_I = -(-1)^r |I|! \, (|I|\!-\!1)! \prod_{j=1}^r \frac{1}{I_j! \,
(I_j\!-\!1)!} \prod_{j=1}^{r-1} \frac{1}{I_j \!+\! I_{j+1}}.
\end{equation}
Here, an empty product has to be interpreted as $1$. Note that
$m_{(N)} = 1$ for all $N \ge 1$. In these terms, the coefficient of
the term
$$
P_{2I} (Q_{2a}), \; |I|+a=4
$$
on the right-hand side of \eqref{Q8-main} is given by
\begin{equation}\label{mult}
-(-1)^a m_{(I,a)},
\end{equation}
and \eqref{Q8-main} can be stated as
\begin{equation}\label{rec-Q8}
\sum_{|I| + a = 4} (-1)^a m_{(I,a)} P_{2I}(Q_{2a}) = 3!4! 2^8 w_8.
\end{equation}
This is a special case of Conjecture 9.2 in \cite{juhl-power}.

Of course, Theorem \ref{main} does not yet provide an explicit
formula for $Q_8$ in terms of the metric. Such a formula can be
derived by combining it with formulas for the lower-order
GJMS-operators $P_2, P_4, P_6$ and the lower-order $Q$-curvatures
$Q_2, Q_4, Q_6$ in dimension $n=8$. The relevant formulas will be
discussed in Section \ref{final}. However, we emphasize that the
resulting identities for $Q_8$ are structurally less natural than
the description \eqref{Q8-main}.

A second feature which distinguishes the formula \eqref{Q8-main} for
$Q_8$ from other formulas is its universality in the dimension of
the underlying space.

\begin{thm}\label{Q8-univ} On any Riemannian manifold of dimension
$n \ge 8$, Branson's $Q$-curvature $Q_8$ is given by the recursive
formula \eqref{Q8-main}.
\end{thm}

For a proof of Theorem \ref{Q8-univ} for the round spheres $S^n$ see
\cite{juhl-power}. The results of \cite{JK} also cover the formula
\eqref{Q8-main} for the conformally flat M\"obius spheres $(S^q
\times S^p, g_{S^q} - g_{S^p})$ with the round metrics on the
factors.

Now we return to the critical case. A closer examination of
\eqref{Q8-main} shows that in the sum on the right-hand side a
substantial number of cancellations takes place. This leads to the
following result.

\begin{thm}\label{reduced} On manifolds of dimension $8$, Branson's
$Q$-curvature $Q_8$ equals the sum of
\begin{multline}\label{Q8-reduced}
-3 P_2^0(Q_6) - 3 P_6^0(Q_2) + 9 P_4^0(Q_4) \\
+ 8 P_2^0 P_4(Q_2) - 12 P_2^0 P_2 (Q_4) + 12 P_4^0 P_2 (Q_2) - 18
P_2^0 P_2^2 (Q_2),
\end{multline}
the divergence term
\begin{equation}\label{div-form}
6 \delta \left( c (2Q_4 + 3P_2 (Q_2),Q_2) \right),
\end{equation}
where $c(f,g) = f d g - g d f \in \Omega^1(M)$, and
$$
3!4!2^7 v_8.
$$
Here $P_{2N}^0$ denotes the non-constant part of $P_{2N}$.
\end{thm}

The reader should note the tiny difference in the last terms in the
formulas in Theorem \ref{main} and Theorem \ref{reduced}: $2w_8$ is
replaced by $v_8$. Note also that $P_2^0 = \Delta$.

Since the operators $P_{2N}^0$ are of the form $\delta( S_{2N} d)$
for some geometric operators $S_{2N}$ on $\Omega^1(M)$ (see
\cite{bran-2}), Theorem \ref{reduced} reproves the following special
case of a result of Graham and Zworski (see \cite{GZ}).

\begin{corr}\label{total} On closed manifolds $M$ of dimension $8$,
\begin{equation}\label{integral}
\int_M Q_8 vol = 3!4! 2^7 \int_M v_8 vol.
\end{equation}
\end{corr}

The present paper rests on the approach to $Q$-curvature developed
in \cite{juhl-book}. For full details we refer to this book and to
Chapter 1 of \cite{BJ}.\footnote{We also use the opportunity to
correct some misprints in \cite{juhl-book}.}

The paper is organized as follows. Section \ref{proof1} contains a
proof of the recursive formula \eqref{Q8-main} for the critical
$Q$-curvature $Q_8$. The key observation in this proof is the
vanishing of the polynomial $\V_8(\lambda)$ (Proposition
\ref{V8-van}). Although the proof of Theorem \ref{main} only applies
the vanishing of the leading coefficient of $\V_8(\lambda)$, the
vanishing result Proposition \ref{V8-van} is of independent
interest. It shows that the vanishing property used here is related
to other relations as e.g. the holographic formula for $Q_8$ (see
\cite{gj}). The reduced form \eqref{Q8-reduced} of \eqref{Q8-main}
is derived in Section \ref{proof2}. In Section \ref{universal}, we
prove the universality of \eqref{Q8-main}, i.e., Theorem
\ref{Q8-univ}. This proof also sheds new light on the proof in
Section \ref{proof1}. In dimension $n \ge 8$, it is still true that
an analog of the polynomial $\V_8(\lambda)$ has a vanishing leading
coefficient although the polynomial itself does not vanish. This
fact can be used to extend the arguments of Section \ref{proof1}.
Here we give an alternative argument. The central point is to prove
the formula in Proposition \ref{eval}. This identity will be
established as a consequence of a more general result (Theorem
\ref{fundamental}) which also provides a certain explanation of the
appearance of the square root of $v(r)$ in Theorems \ref{main} and
\ref{Q8-univ}. In Section \ref{final}, we discuss analogous
descriptions of $Q_4$ and $Q_6$, and display universal formulas for
the GJMS-operators $P_4$ and $P_6$. Section \ref{final-co} contains
comments on further developments. In particular, we describe the
status of Conjecture 9.2 of \cite{juhl-power}.

\section{Proof of Theorem \ref{main}}\label{proof1}

The basic idea of the proof of Theorem \ref{main} is to compare two
different evaluations of the leading coefficient of the
$Q$-curvature polynomial $Q_8^{res}(\lambda)$.

We first recall the notion of $Q$-curvature polynomials (or
$Q$-polynomials for short) as introduced in \cite{juhl-book}; see also
Section 1.6 of \cite{BJ}. For this purpose, let
\begin{equation}\label{expansion}
u \sim \sum_{N \ge 0} r^{\lambda+2N} \T_{2N}(g;\lambda)(f), \;
\T_0(g;\lambda)(f) = f, \; r \to 0
\end{equation}
be the asymptotic expansions of an eigenfunction $u$ of the
Laplace-Beltrami operator of a Poincar\'e-Einstein metric $g_+$ as
in \eqref{PE1} and \eqref{PE2}:
$$
-\Delta_{g_+} u = \lambda(n-\lambda)u.
$$
In the asymptotic expansion \eqref{expansion}, we suppresses the
analogous contributions of the form $\sum_{N \ge 0} r^{n-\lambda+2N}
b_{2N}$. The coefficients $\T_{2N}(g;\lambda)$ are rational families
(in $\lambda$) of differential operators of the form
\begin{equation}\label{TP}
\T_{2N}(g;\lambda) = \frac{1}{2^{2N} N! (\f\!-\!\lambda\!-\!1)
\cdots (\f\!-\!\lambda\!-\!N)} P_{2N}(g;\lambda)
\end{equation}
with a polynomial family $P_{2N}(g;\lambda) = \Delta^N + LOT$. In
particular, the poles of $\T_{2N}(\lambda)$ are contained in the set
$\left\{\f-1,\dots,\f-N\right\}$. The families $P_{2N}(g;\lambda)$
contain the GJMS-operators for special parameters $\lambda$. More
precisely, we have (see \cite{GZ})
\begin{equation}\label{spectral}
P_{n-2N}(g;N) = P_{n-2N}(g) \quad \mbox{for $N=0,1,\dots,\f$.}
\end{equation}

\begin{defn}[{\bf $Q$-curvature polynomials}] For even $n \ge 2$
and $2 \le 2N \le n$, the $N^{th}$ $Q$-curvature polynomial is
defined by
\begin{multline}\label{Q-pol}
Q_{2N}^{res}(g;\lambda) = -2^{2N} N!
\left(\left(\lambda\!+\!\f\!-\!2N\!+\!1\right) \cdots
\left(\lambda\!+\!\f\!-\!N\right)\right)
\\ \times \left[\T_{2N}^*(g;\lambda\!+\!n\!-\!2N)(v_0) + \dots +
\T_0^*(g;\lambda\!+\!n\!-\!2N)(v_{2N})\right].
\end{multline}
We also set $Q_0^{res}(g;\lambda)=-1$.
\end{defn}

In particular, the critical $Q$-curvature polynomial is given by the
formula
\begin{multline}\label{Q-pol-crit}
Q_n^{res}(g;\lambda) = - 2^n \left(\f\right)!
\left(\left(\lambda\!-\!\f\!+\!1\right) \cdots \lambda\right)
\left[\T_n^*(g;\lambda)(v_0) + \dots +
\T_0^*(g;\lambda)(v_n)\right].
\end{multline}

$Q_{2N}^{res}(\lambda)$ is a polynomial of degree $N$. In
particular, the critical $Q$-curvature polynomial
$Q_n^{res}(\lambda)$ has degree $\f$. It has vanishing constant
term, i.e., $Q_n^{res}(0)=0$, and satisfies
\begin{equation}\label{name}
\dot{Q}_n^{res}(0) = Q_n.
\end{equation}
It is the latter property which motivates the name. For proofs of
\eqref{name} see \cite{gj}, \cite{BJ} or \cite{juhl-book}.

The polynomials $Q_{2N}^{res}(\lambda)$ are proportional to the
constant terms of the so-called residue families
$D_{2N}^{res}(\lambda)$:
$$
Q_{2N}^{res}(\lambda) = -(-1)^N D_{2N}^{res}(\lambda)(1).
$$
One of the basic properties of the families $D_{2N}^{res}(\lambda)$
is that, for special values of the parameter $\lambda$, they factor
into products of lower order residue families and GJMS-operators.
These factorization identities allow to express any of these
families in terms of respective lower order families and
GJMS-operators. As a consequence, we have

\begin{prop}\label{q-factor} The $Q$-curvature polynomials
$Q_{2N}^{res}(\lambda)$ satisfy the factorization relations
$$
Q_{2N}^{res}\left(-\f+2N-j\right) = (-1)^j P_{2j} \left(
Q_{2N-2j}^{res}\left(-\f+2N-j\right)\right)
$$
for $j=1,\dots,N$.
\end{prop}

In particular, the critical $Q$-curvature polynomial
$Q_8^{res}(\lambda)$ satisfies the relations
\begin{align*}
Q_8^{res}(3) & = - P_2 \left(Q_6^{res}(3)\right), \\
Q_8^{res}(2) & = P_4 \left(Q_4^{res}(2)\right), \\
Q_8^{res}(1) & = - P_6 \left(Q_2^{res}(1)\right), \\
Q_8^{res}(0) & = -P_8 (1) = 0.
\end{align*}
Since $Q_8^{res}(\lambda)$ is a polynomial of degree $4$, these
relations together with
\begin{equation}\label{dot8}
\dot{Q}^{res}_8(0) = Q_8
\end{equation}
imply that the coefficients of $Q_8^{res}(\lambda)$ can be written
as linear combinations of
$$
P_2 \left(Q_6^{res}(3)\right), \; P_4 \left(Q_4^{res}(2)\right), \;
P_6 \left(Q_2^{res}(1)\right) \quad \mbox{and} \quad  Q_8.
$$
Now we have $Q_2^{res}(\lambda)=\lambda Q_2$, and Theorem 6.11.8 in
\cite{juhl-book} yields the formulas
\begin{equation}\label{Q4-pol}
Q_4^{res}(\lambda) = -\lambda(\lambda+1) Q_4 - \lambda(\lambda+2)
P_2(Q_2)
\end{equation}
and
\begin{equation}\label{Q6-pol}
Q_6^{res}(\lambda) = \frac{1}{2} \lambda^2 (\lambda\!-\!1) Q_6 +
\lambda^2 (\lambda\!+\!1) P_2 \left(Q_4 + \frac{3}{2}
P_2(Q_2)\right) - \lambda(\lambda\!+\!1)(\lambda\!-\!1) P_4(Q_2).
\end{equation}

It is easy to verify that $Q_4^{res}(\lambda)$ and
$Q_6^{res}(\lambda)$ indeed satisfy the respective relations
\begin{align*}
Q_6^{res}(1) & = - P_2 \left(Q_4^{res}(1)\right), \\
Q_6^{res}(0) & = P_4 \left(Q_2^{res}(0)\right), \\
Q_6^{res}(-1) & = P_6 (1) = -Q_6
\end{align*}
and
\begin{align*}
Q_4^{res}(-1) = - P_2 \left(Q_2^{res}(-1)\right), \\
Q_4^{res}(-2) = - P_4 (1) = -2 Q_4.
\end{align*}

Here some comments are in order. The above $3$ factorization
relations for the cubic polynomial $Q^{res}_6(\lambda)$ do not
suffice for its characterization. However, in dimension $n \ge 10$,
the analogous relations together with $Q_6^{res}(0)=0$ yield a
characterization, and the above formula for $n=8$ follows by
``analytic continuation'' in the dimension. This argument differs
from that in \cite{juhl-book}, where the proof of the analog of
\eqref{Q6-pol} for general dimensions rests on an explicit formula
for $Q_6$, and the vanishing property $Q_6^{res}(0)=0$ appears as a
consequence. Similarly, the two relations for the quadratic
polynomial $Q_4^{res}(\lambda)$ together with the vanishing property
$Q_4^{res}(0)=0$ characterize this polynomial. For a proof of
$Q_{2N}^{res}(0)=0$ in full generality see Section 1.6 of \cite{BJ}.

Now the above formulas for $Q_2^{res}(\lambda)$,
$Q_4^{res}(\lambda)$ and $Q_6^{res}(\lambda)$ together with
\eqref{dot8} yield a formula for the quartic polynomial
$Q_8^{res}(\lambda)$. The actual calculation shows that the
coefficient of $\lambda^4$ in $-6 Q_8^{res}(\lambda)$ is given by
\begin{multline}\label{result-1}
Q_8 + 3 P_2(Q_6) + 3 P_6(Q_2) - 9 P_4(Q_4) \\
- 8 P_2 P_4(Q_2) + 12 P_2^2 (Q_4) - 12 P_4 P_2 (Q_2) + 18
P_2^3(Q_2).
\end{multline}

We compare formula \eqref{result-1} with the result of a direct
evaluation of the leading coefficient of $-6Q_8^{res}(\lambda)$. In
this direction, we first observe that the definitions imply
\begin{multline}\label{first-Q8}
-Q_8^{res}(\lambda) = P_8^*(\lambda)(1) - 16 \lambda
P_6^*(\lambda)(v_2) + 2^6 3 \lambda(\lambda\!-\!1)
P_4^*(\lambda)(v_4) \\ - 2^9 3 \lambda(\lambda\!-\!1)(\lambda\!-\!2)
P_2^*(\lambda)(v_6) + 2^{10} 3!
\lambda(\lambda\!-\!1)(\lambda\!-\!2)(\lambda\!-\!3) v_8.
\end{multline}
We avoid a consideration of the complicated family $P_8(\lambda)$ by
showing that the term $P_8^*(\lambda)(1)$ can be expressed as a
linear combination of the other four terms in \eqref{first-Q8}. For
this purpose, we define for any manifold of even dimension $n$ the
polynomial
\begin{equation}\label{V-general}
\V_n(\lambda) = \left[\lambda (\lambda\!-\!1) \cdots
\left(\lambda\!-\!\f\!+\!1\right)\right] \sum_{j=0}^\f (n\!+\!2j)
\T_{2j}^*(\lambda) (v_{n-2j}).
\end{equation}
The polynomial $\V_n$ has degree $\f$. In \cite{juhl-book}, we
formulated the conjecture that $\V_n(\lambda)$ vanishes (see
Conjecture 6.11.2). The vanishing of $\V_n(\lambda)$ is equivalent
to the system of $\f+1$ relations which expresses the vanishing of
its coefficients. Here we are interested in the quartic polynomial
$\V_8(\lambda)$. Only the vanishing of its leading coefficient will
be important in the sequel. We establish this vanishing as a
consequence of the vanishing of $\V_8(\lambda)$. The proof of this
property will {\em not} require to make explicit the family
$P_8(\lambda)$.

For an alternative proof of the vanishing of the leading coefficient
of an analog of $\V_8(\lambda)$ for manifolds of general dimensions
we refer to Section \ref{universal}.

The following result confirms Conjecture 6.11.2 of \cite{juhl-book}
for $\V_8(\lambda)$.

\begin{prop}\label{V8-van} For any manifold of dimension $n=8$, the
quartic polynomial
\begin{equation}\label{V8}
\V_8(\lambda) =
\left[\lambda(\lambda\!-\!1)(\lambda\!-\!2)(\lambda\!-\!3)\right]
\sum_{j=0}^4 (8+2j) \T_{2j}^*(\lambda)(v_{8-2j})
\end{equation}
vanishes identically.
\end{prop}

The strategy of the proof of Proposition \ref{V8-van} is the
following. For general even $n$, the vanishing of $\V_n(\lambda)$ is
equivalent to the conditions
\begin{equation}\label{van-cond}
\V_n(0) = \V_n(1) = \dots = \V_n\left(\f-1\right) = 0 \quad
\mbox{and} \quad \dot{\V}_n(0) = 0.
\end{equation}
But the conditions
\begin{equation}\label{two-cond}
\V_n(0)=0 \quad \mbox{and} \quad \dot{\V}_n(0)=0
\end{equation}
are known to be satisfied in full generality (Theorem 6.11.12 in
\cite{juhl-book}). We prove that the remaining conditions follow
from a simpler system of conditions, and verify the latter ones for
$\V_8(\lambda)$.

For the convenience of the reader, we also describe the arguments
which prove \eqref{two-cond}. First, \eqref{TP} shows that
$\V_n(0)=0$ is equivalent to $P_n^*(0)(1) = 0$. Thus, the first
assertion in \eqref{two-cond} follows from $P_n(0) = P_n$ and the
fact that $P_n$ is a self-adjoint operator with vanishing constant
term. The second condition in \eqref{two-cond} is more subtle.
\eqref{TP} shows that the linear coefficient of $\V_n(\lambda)$ is
given by
\begin{equation}\label{linear}
\frac{(-1)^\f}{2^{n-1}(\f)!} \left( n \dot{P}_n^*(0)(v_0) - 2^{n-1}
\left(\f\!-\!1\right)! \left(\f\right)! \sum_{j=0}^{\f-1} (n\!+\!2j)
\T_{2j}^*(0)(v_{n-2j})\right).
\end{equation}
Now we combine $\dot{P}_n (0)(1) = (-1)^\f Q_n$ (\cite{GZ}), the
relation
$$
n \dot{P}_n^*(0)(1) = n \dot{P}_n(0)(1) + 2^n \left(\f\right)!
\left(\f\!-\!1\right)! \sum_{j=0}^{\f-1} 2j \T_{2j}^*(0) (v_{n-2j})
$$
(see \cite{gj}, Proposition 2) and the {\em holographic formula}
\begin{equation}\label{holo}
n (-1)^\f Q_n = 2^{n-1} \left(\f\right)! \left(\f\!-\!1\right)!
\sum_{j=0}^{\f-1} (n\!-\!2j) \T_{2j}^*(0) (v_{n-2j})
\end{equation}
for $Q_n$ (see \cite{gj}, Theorem 1). It follows that
\begin{multline*}
n \dot{P}_n^*(0)(1) = 2^{n-1} \left(\f\right)!
\left(\f\!-\!1\right)! \sum_{j=0}^{\f-1} (n\!-\!2j) \T_{2j}^*(0) (v_{n-2j})
\\[-2mm]
+ 2^n \left(\f\right)! \left(\f\!-\!1\right)! \sum_{j=0}^{\f-1} 2j
\T_{2j}^*(0) (v_{n-2j}),
\end{multline*}
i.e., \eqref{linear} vanishes. This proves the second assertion.

The following result provides a sufficient condition for the
remaining vanishing properties in \eqref{van-cond}.

\begin{lemm}\label{van-equiv} For $N=1,\dots,\f-1$, the condition
$$
\V_n(N)=0
$$
follows from
$$
\sum_{j=0}^N (n\!-\!j) \T_{2N-2j}^*(n\!-\!N)(v_{2j}) = 0.
$$
\end{lemm}

\begin{proof} By \eqref{TP}, the vanishing of $\V_n(N)$ is
equivalent to
\begin{multline}\label{A}
2n \frac{(-1)^\f}{2^n \left(\f\right)!} P_n^*(N)(v_0) + (2n\!-\!2)
\frac{(-1)^{\f-1}}{2^{n-2} \left(\f\!-\!1\right)!} N
P_{n-2}^*(N)(v_2) \\ + \cdots + (2n\!-\!2N)
\frac{(-1)^{\f-N}}{2^{n-2N} \left(\f\!-\!N\right)!} N!
P_{n-2N}^*(N)(v_{2N}) = 0.
\end{multline}
Now we have the factorization relations
\begin{equation}\label{factor}
P_{n-2j}(N) = P_{2N-2j}(n\!-\!N) P_{n-2N}, \; j=0,\dots,N.
\end{equation}
These follow from \eqref{spectral} and the fact that, for certain
values of $\lambda$, the operators $\T_{2N}(\lambda)$ can be
determined in stages, i.e., the map
$$
f r^N \mapsto \T_{n-2N}(N)(f) r^{n-N} \mapsto \T_{2N-2j}(n\!-\!N)
\left( \T_{n-2N}(N)(f)\right) r^{n-2j+N}
$$
coincides with
$$
f r^N \mapsto \T_{n-2j}(N)(f) r^{n-2j+N}
$$
Eq.~\eqref{factor} show that \eqref{A} follows from
\begin{multline}\label{B}
2n \frac{(-1)^\f}{2^n \left(\f\right)!} P_{2N}^*(n\!-\!N)(v_0) +
(2n\!-\!2)  \frac{(-1)^{\f-1}}{2^{n-2} \left(\f\!-\!1\right)!} N
P_{2N-2}^*(n\!-\!N)(v_2) \\ + \cdots + (2n\!-\!2N)
\frac{(-1)^{\f-N}}{2^{n-2N} \left(\f\!-\!N\right)!} N! v_{2N} = 0.
\end{multline}
But \eqref{TP} implies that this relation is equivalent to
$$
2^{2N-n} \frac{N!}{(\f\!-\!N)!} \left( 2n \T_{2N}^*(n\!-\!N)(v_0) +
\cdots + (2n\!-\!2N) \T_0^*(n\!-\!N)(v_{2N}) \right) = 0.
$$
The proof is complete.
\end{proof}

Now we are ready to prove Proposition \ref{V8-van}.

\begin{proof}[Proof of Proposition \ref{V8-van}] It only remains to prove
that $\V_8(1)=\V_8(2)=\V_8(3)=0$. In order to reduce the amount of
numerical factors, we actually prove the more general result that
$$
\V_n(1) = \V_n(2) = \V_n(3) = 0
$$
for all $n \ge 8$.\footnote{The arguments also prove that
$\V_n(1)=\V_n(2)=0$ for $n \ge 6$.} By Lemma \ref{van-equiv}, these
relations follow from the identities
\begin{equation}\label{a}
n \T_2^*(n\!-\!1)(v_0) + (n\!-\!1) v_2 = 0,
\end{equation}
\begin{equation}\label{b}
n \T_4^*(n\!-\!2)(v_0) + (n\!-\!1) \T_2^*(n\!-\!2)(v_2) + (n\!-\!2)
v_4 = 0
\end{equation}
and
\begin{equation}\label{c}
n \T_6^*(n\!-\!3)(v_0) + (n\!-\!1) \T_4^*(n\!-\!3)(v_2) + (n\!-\!2)
\T_2^*(n\!-\!3)(v_4) + (n\!-\!3) v_6 = 0.
\end{equation}
In the remainder of the proof we confirm these three relations. For
this purpose, we apply the following explicit formulas for the
quantities involved. First of all, the families $P_2(\lambda)$ and
$P_4(\lambda)$ are given by
\begin{equation}\label{P2}
P_2(\lambda) = \Delta - \lambda \J
\end{equation}
and
\begin{multline}\label{P4}
P_4(\lambda) = (\Delta \!-\! (\lambda\!+\!2)\J)(\Delta \!-\! \lambda
\J) + \lambda(2\lambda\!-\!n\!+\!2) |\Rho|^2 \\ +
2(2\lambda\!-\!n\!+\!2) \delta (\Rho d) + (2\lambda\!-\!n\!+\!2)
(d\J,d).
\end{multline}
These two formulas are contained in Theorem 6.9.4 in
\cite{juhl-book}. Next, we have
\begin{equation}\label{v24}
v_2 = -\frac{1}{2} \J \quad \mbox{and} \quad v_4 = \frac{1}{8}
(\J^2-|\Rho|^2).
\end{equation}
Eq.~\eqref{a} easily follows from the definitions. Eq.~\eqref{b} is
equivalent to
\begin{multline*}
\frac{n}{8(n\!-\!2)n} \left( -n(\Delta\!-\!(n\!-\!2)\J) \J +
(n\!-\!2)^2 |\Rho|^2 - (n\!-\!2) \Delta \J \right) \\
+ \frac{1}{n\!-\!2} (\Delta \J\!-\!(n\!-\!2) \J^2) + \frac{n\!-\!2}{8}
(\J^2\!-\!|\Rho|^2) = 0.
\end{multline*}
It is straightforward to verify this relation. The proof of
\eqref{c} requires some more work. We start by observing that
\eqref{c} is equivalent to
\begin{multline}\label{d}
P_6^*(n\!-\!3)(v_0) - 6 (n\!-\!1) P_4^*(n\!-\!3)(v_2) \\ + 24 (n\!-\!2)^2
P_2^*(n\!-\!3)(v_4) - 48 (n\!-\!2)(n\!-\!3)(n\!-\!4) v_6 = 0.
\end{multline}
For the evaluation of \eqref{d}, we apply a formula for
$P_6(\lambda)$ which was derived as formula (6.10.2) in
\cite{juhl-book}. Its formulation requires to introduce some
notation. Let
$$
g_t = g + t g_{(2)} + t^2 g_{(4)} = g - t \Rho + \frac{1}{4} t^2
\left( \Rho^2 - \frac{\B}{n-4} \right),
$$
where $\B$ is the Bach tensor (see \eqref{bach}). Iterated
derivatives with respect to $t$ (at $t=0$) will be denoted by $'$.
In particular, $\Delta'$ and $\Delta''$ are the first and second
metric variations of $\Delta$ for the variation of $g$ defined by
$g_t$. In these terms,
\begin{multline}\label{p6-form}
P_6(\lambda) u = 4 (n\!-\!4\!-\!2\lambda)(n\!-\!2\!-\!2\lambda)
\left[\lambda (\log \det g)''' + \Delta'' \right] (u) \\
+ 4 (n\!-\!4\!-\!2\lambda) \left[ (\lambda\!+\!2) (\log \det g)'' +
\Delta' \right] P_2(\lambda) (u) \\
+ (\Delta \!-\! (\lambda\!+\!4)\J) P_4(\lambda) u.
\end{multline}
For more details see Section \ref{universal}. Eq.~\eqref{p6-form}
implies
\begin{multline}\label{e}
P_6^*(\lambda)(1) = 4 (n\!-\!4\!-\!2\lambda)(n\!-\!2\!-\!2\lambda)
\left[\lambda (\log \det g)''' + {\Delta''}^*(1) \right] \\
+ 4 (n\!-\!4\!-\!2\lambda) P_2^*(\lambda) \left[ (\lambda\!+\!2)
(\log \det g)'' + \Delta'^* (1) \right] \\
+ P_4^*(\lambda) (\Delta \!-\! (\lambda\!+\!4)\J)(1).
\end{multline}
In order to determine the quantity ${\Delta''}^*(1)$, we combine
\eqref{e} with the relation
$$
P_6\left(\f\!-\!3\right) = P_6
$$
and the fact that $P_6$ is self-adjoint. This yields the formula
\begin{equation}
4 \left( \Delta'' - (\Delta'')^* \right) (1) = - \Delta |\Rho|^2 + 4
\delta (\Rho d\J).
\end{equation}
For $n=6$, the details of the calculation can be found in Section
6.10 of  \cite{juhl-book}. As to be expected, the result does not
depend on the dimension. Since $\Delta''(1) = 0$, we find
\begin{equation}\label{d2}
4 (\Delta'')^* (1) = \Delta |\Rho|^2 - 4 \delta (\Rho d\J).
\end{equation}
Now, we evaluate \eqref{e} by using \eqref{d2},
\begin{align}
(\log \det g)'' & = -\frac{1}{2} |\Rho|^2, \label{second}\\
(\log \det g)'''& = -\frac{1}{2(n\!-\!4)} (\B,\Rho) - \frac{1}{2}
\tr (\Rho^3), \label{third}
\end{align}
and
$$
(\Delta')^*(1) = \frac{1}{2} \Delta \J.
$$
For the proofs of these results we refer to Section 6.10 of
\cite{juhl-book}. We obtain
\begin{multline*}
P_6^*(\lambda)(1) = (n\!-\!4\!-\!2\lambda)(n\!-\!2\!-\!2\lambda)
\left[ \Delta |\Rho|^2 - 4 \delta (\Rho d\J) - \frac{2 \lambda}{n-4}
(\B,\Rho) -2\lambda \tr(\Rho^3) \right] \\
+ 2 (n\!-\!4\!-\!2\lambda) \left[ -(\lambda\!+\!2)(\Delta-\lambda
\J) |\Rho|^2 + (\Delta\!-\!\lambda\J) \Delta \J \right]
-(\lambda\!+\!4) P_4^*(\lambda) \J.
\end{multline*}
In particular, we find
\begin{multline*}
P_6^*(n\!-\!3)(1) = (n\!-\!2)(n\!-\!4) \left[ \Delta |\Rho|^2 -4
\delta(\Rho d\J)
-\frac{2(n\!-\!3)}{n\!-\!4} (\B,\Rho) - 2(n\!-\!3) \tr(\Rho^3) \right] \\
-2(n\!-\!2) \left[ -(n\!-\!1) (\Delta \!-\! (n\!-\!3) \J) |\Rho|^2 +
(\Delta\!-\!(n\!-\!3)\J) \Delta \J \right] - (n\!+\!1) P_4^*(n\!-\!3) \J.
\end{multline*}
Thus, the left-hand side of \eqref{d} equals the sum of
\begin{multline*}
(n\!-\!2)(n\!-\!4) \left[ \Delta |\Rho|^2 -4 \delta(\Rho d\J)
-\frac{2(n\!-\!3)}{n\!-\!4} (\B,\Rho) - 2(n\!-\!3) \tr(\Rho^3) \right] \\
-2(n\!-\!2) \left[ -(n\!-\!1) (\Delta \!-\! (n\!-\!3) \J) |\Rho|^2 +
(\Delta\!-\!(n\!-\!3)\J) \Delta \J \right] + 2(n\!-\!2)  P_4^*(n\!-\!3) \J,
\end{multline*}
where
\begin{multline*}
P_4^*(n\!-\!3)\J = (\Delta\!-\!(n\!-\!3)\J)(\Delta\!-\!(n\!-\!1)\J)\J \\
+ (n\!-\!3)(n\!-\!4) |\Rho|^2 \J + 2(n\!-\!4) \delta(\Rho d\J) +
(n\!-\!4) \delta(\J d\J),
\end{multline*}
and
\begin{multline*}
3 (n\!-\!2)^2 (\Delta\!-\!(n\!-\!3)\J) (\J^2 \!-\! |\Rho|^2) \\
-48 (n\!-\!2)(n\!-\!3)(n\!-\!4) \left( -\frac{1}{8} \tr (\wedge^3 \Rho) -
\frac{1}{24(n\!-\!4)} (\B,\Rho) \right).
\end{multline*}
Here we made use of formula \eqref{v6}. Now a direct calculation
using Newton's formula
\begin{equation}\label{newton}
6 \tr(\wedge^3 \Rho) = \J^3 - 3 \J |\Rho|^2 + 2 \tr (\Rho^3)
\end{equation}
shows that this sum vanishes.
\end{proof}

We continue with the

\begin{proof}[Proof of Theorem \ref{main}] The vanishing of
$\V_8(\lambda)$ is equivalent to the identity\footnote{The
identities \eqref{first-Q8} and \eqref{first-P8} can be found also
in the proof of Theorem 6.13.1 in \cite{juhl-book}. Here we correct
misprints in the coefficients of $v_8$ in both formulas.}
\begin{multline}\label{first-P8}
P_8^*(\lambda)(1) = 14 \lambda P_6^*(\lambda)(v_2) - 2^4 9
\lambda(\lambda\!-\!1) P_4^*(\lambda)(v_4) \\ + 2^6 15
\lambda(\lambda\!-\!1)(\lambda\!-\!2) P_2^*(\lambda)(v_6) - 2^{10} 3
\lambda(\lambda\!-\!1)(\lambda\!-\!2)(\lambda\!-\!3) v_8.
\end{multline}
Combining this result with \eqref{first-Q8} yields
\begin{multline*}
-Q_8^{res}(\lambda) = - 2\lambda P_6^*(\lambda)(v_2) + 2^4 3
\lambda(\lambda\!-\!1) P_4^*(\lambda)(v_4) \\ - 2^6 3^2
\lambda(\lambda\!-\!1)(\lambda\!-\!2) P_2^*(\lambda) (v_6) + 2^{10}
3 \lambda(\lambda\!-\!1)(\lambda\!-\!2)(\lambda\!-\!3)v_8.
\end{multline*}
From this formula, we read-off the coefficient of $\lambda^4$ in $-6
Q_8^{res}(\lambda)$ as
\begin{equation}\label{intermediate}
-12 P_6^*(\lambda)^{[3]}(v_2) + 2^5 3^2 P_4^*(\lambda)^{[2]}(v_4) -
2^7 3^3 P_2^*(\lambda)^{[1]} (v_6) + 2^{11} 3^2 v_8,
\end{equation}
where the superscripts indicate the coefficients of the respective
powers of $\lambda$. Now we have
$$
P_2^*(\lambda)^{[1]} = - \J = 2v_2, \quad P_4^*(\lambda)^{[2]} =
\J^2 + 2 |\Rho|^2 = -16 v_4 + 12 v_2^2
$$
by \eqref{P2} and \eqref{P4}, and
$$
P_6^*(\lambda)^{[3]} = 16 (\log \det g)''' + 8 (\log \det g)'' \J -
\J (\J^2+2|\Rho|^2)
$$
by \eqref{p6-form}. A calculation using \eqref{second},
\eqref{third}, \eqref{v24} and \eqref{v6} shows that
\begin{align*}
(\log \det g)'' & =  4v_4 - 2 v_2^2, \\
(\log \det g)''' & = 12 v_6 - 12 v_2v_4 + 4v_2^3.
\end{align*}
Hence we find
\begin{equation}\label{P6-leading}
P_6^*(\lambda)^{[3]}(v_2) = 24 (8v_2v_6 - 12 v_2^2 v_4 + 5 v_2^4),
\end{equation}
and it follows that \eqref{intermediate} is given by the sum of
\begin{equation}\label{start-h}
2^5 3^2 \left( -32 v_2 v_6 + 24 v_2^2 v_4 - 5 v_2^4 - 16 v_4^2
\right)
\end{equation}
and $3! 4! 2^7 v_8$. Now Lemma \ref{sroot} shows that the
coefficient of $\lambda^4$ in $-6 Q_8^{res}(\lambda)$ equals
$$
2^{12} 3^2 w_8 = 3! 4! 2^8 w_8.
$$
Comparing this result with \eqref{result-1} implies the assertion.
\end{proof}

The following elementary algebraic result was used in the proof of
Theorem \ref{main}.

\begin{lemm}\label{sroot} Let
$$
1 + w_2 r^2 + w_4 r^4 + w_6 r^6 + w_8 r^8 + \cdots
$$
be the Taylor series of the function $w(r) = \sqrt{v(r)}$ with
$$
v(r) = 1 + v_2 r^2 + v_4 r^4 + v_6 r^6 + v_8 r^8 + \cdots.
$$
Then
\begin{align*}
2 w_2 & = v_2, \\
2 w_4 & = \frac{1}{4} (4 v_4 - v_2^2), \\
2 w_6 & = \frac{1}{8} (8 v_6 - 4 v_4 v_2 + v_2^3), \\
2 w_8 & = \frac{1}{64} (64 v_8 - 32 v_6v_2 - 16 v_4^2 + 24 v_2^2 v_4
- 5 v_2^4).
\end{align*}
\end{lemm}

The assertion follows by squaring the Taylor series of $w$.

\section{Proof of Theorem \ref{reduced}}\label{proof2}

We derive Theorem \ref{reduced} from Theorem \ref{main}. The proof
consists of two steps. In the first step we establish

\begin{prop}\label{Qv} On manifolds of dimension $n \ge 8$,
\begin{multline}\label{vw}
3! 4! 2^8 w_8 - 3!4!2^7 v_8 \\ = -12 \left[ Q_6 + 2P_2(Q_4) -
2P_4(Q_2) + 3 P_2^2(Q_2) \right] Q_2 - 18 \left[Q_4 +
P_2(Q_2)\right]^2.
\end{multline}
\end{prop}

\begin{proof} By \eqref{v8w8}, the assertion is equivalent to
\begin{multline*}
48 (32 v_6 v_2 + 16 v_4^2 - 24 v_2^2 v_4 + 5v_2^4) \\
= 2 \left[Q_6 + 2P_2(Q_4) - 2P_4(Q_2) + 3P_2^2(Q_2) \right] Q_2 + 3
\left[ Q_4 + P_2(Q_2) \right]^2.
\end{multline*}
For the proof of this identity we regard the relations
$$
Q_2 = -2 v_2, \quad Q_4 = -P_2(Q_2) - Q_2^2 + 16 v_4
$$
and
$$
Q_6 = \left[-2P_2(Q_4) + 2P_4(Q_2) - 3P_2^2(Q_2)\right] - 6 \left[
Q_4 + P_2(Q_2) \right] Q_2 - 2! 3! 2^5 v_6
$$
in dimension $n \ge 8$ (Proposition \ref{Q4-both} and Proposition
\ref{Q6-both}) as formulas for $v_2$, $v_4$ and $v_6$, and find
\begin{multline*}
48 (32 v_6 v_2 + 16 v_4^2 - 24 v_2^2 v_4 + 5v_2^4) \\
= 2 \left[Q_6 \!+\! 2P_2(Q_4) \!-\! 2P_4(Q_2) \!+\!
3P_2^2(Q_2)\right] Q_2 + 12 \left[Q_4\!+\!P_2(Q_2)\right] Q_2^2  \\
+ 3 \left[ Q_4\!+\!P_2(Q_2)\!+\!Q_2^2 \right]^2 - 18 \left[
Q_4\!+\!P_2(Q_2)\!+\!Q_2^2 \right] Q_2^2 + 15 Q_2^4.
\end{multline*}
From here the assertion follows by simplification.
\end{proof}

We emphasize the important structural fact that the terms
$$
Q_6 + 2P_2(Q_4) - 2P_4(Q_2) + 3P_2^2(Q_2) \quad \mbox{and} \quad Q_4
+ P_2(Q_2)
$$
on the right-hand side of \eqref{vw} naturally appear also in the
respective recursive formulas \eqref{rec-Q4a} and \eqref{rec-Q6a}
for $Q_4$ and $Q_6$.

Now, Theorem \ref{main} and Proposition \ref{Qv} (in dimension
$n=8$) show that
\begin{multline*}
Q_8 = \big[\!-3 P_2(Q_6) \!-\! 3 P_6(Q_2) \!+\! 9 P_4(Q_4) \\
\!+\! 8 P_2 P_4(Q_2) \!-\! 12 P_2^2 (Q_4) \!+\! 12 P_4 P_2 (Q_2)
\!-\! 18 P_2^3(Q_2)\big] \\ - 12 \left[Q_6 \!+\! 2P_2(Q_4) \!-\!
2P_4(Q_2) \!+\! 3P_2^2(Q_2)\right] Q_2 - 18
\left[Q_4\!+\!P_2(Q_2)\right]^2 + 3!4!2^7 v_8.
\end{multline*}
The decompositions
\begin{equation}\label{deco}
P_2 = P_2^0 - 3 Q_2, \quad P_4 = P_4^0 + 2 Q_4 \quad \mbox{and}
\quad P_6 = P_6^0 - Q_6
\end{equation}
imply that this sum differs from
\begin{multline*}
\big[\!-3 P^0_2(Q_6) - 3 P^0_6(Q_2) +9 P^0_4(Q_4) \\
+ 8 P_2^0 P_4(Q_2) - 12 P_2^0 P_2 (Q_4) + 12 P_4^0 P_2 (Q_2) - 18
P_2^0 P_2^2(Q_2)\big] + 3!4!2^7 v_8
\end{multline*}
by
\begin{multline*}
\big[9 Q_2 Q_6 + 3 Q_6 Q_2 + 18 Q_4^2 \\ - 24 Q_2 P_4(Q_2) + 36 Q_2
P_2(Q_4)+ 24 Q_4 P_2(Q_2) + 54 Q_2 P_2^2(Q_2)\big] \\
- 12 \left[Q_6 \!+\! 2P_2(Q_4) \!-\! 2P_4(Q_2) \!+\!
3P_2^2(Q_2)\right] Q_2 - 18 \left[Q_4\!+\!P_2(Q_2)\right]^2.
\end{multline*}
In the latter sum, we replace $P_2$ and $P_4$ by the decompositions
in \eqref{deco} and simplify. We find
\begin{multline*}
12 \left[ \Delta(Q_4) Q_2 - Q_4 \Delta(Q_2) \right] + 18 \left[
\Delta^2(Q_2) Q_2 - \Delta(Q_2) \Delta(Q_2) \right] \\ + 54 \left[
\Delta(Q_2) Q_2^2 - Q_2 \Delta(Q_2^2) \right].
\end{multline*}
Further simplification yields the formula in Theorem \ref{reduced}.

\section{Proof of Theorem \ref{Q8-univ}}\label{universal}

Let $n \ge 8$. Then the $Q$-curvature polynomial
$Q_8^{res}(\lambda)$ can be written in the form
\begin{multline}\label{Q8-gen}
Q_{8}^{res}(\lambda) = - \lambda \prod_{k=1}^3
\left( \frac{\lambda+\f-8+k}{k} \right) Q_8 \\
+ \lambda \sum_{j=1}^3 (-1)^j \prod_{k=1 \atop k \ne j}^4 \left(
\frac{\lambda+\f-8+k}{k-j} \right) P_{2j} \left(
\QR_{8-2j}^{res}\left(-\f\!+\!8\!-\!j\right)\right),
\end{multline}
where the polynomials $\QR_{2j}^{res}(\lambda)$ are determined by
$$
\lambda \QR^{res}_{2j}(\lambda) = Q_{2j}^{res}(\lambda), \; j=1,2,3.
$$
We recall that $\QR_{2j}^{res}(\lambda)$ is well-defined since by
Theorem 1.6.6 in \cite{BJ}
$$
Q_{2j}^{res}(0)=0, \; j=1,\dots,4.
$$
The formula in \eqref{Q8-gen} follows from the fact that the $4$
factorization identities in Proposition \ref{q-factor} (for $N=4$)
and the vanishing property $Q_8^{res}(0)=0$ characterize this
polynomial. Combining \eqref{Q8-gen} with the analogous formulas
\begin{multline}\label{Q6-gen}
\QR_6^{res}(\lambda) = \prod_{k=1}^2 \left( \frac{\lambda+\f-6+k}{k}
\right) Q_6 \\ + \sum_{j=1}^2 (-1)^j \prod_{k=1 \atop k \ne j}^3
\left( \frac{\lambda+\f-6+k}{k-j} \right) P_{2j} \left(
\QR_{6-2j}^{res}\left(-\f\!+\!6\!-\!j\right)\right)
\end{multline}
and
\begin{equation}\label{Q4-gen}
\QR_4^{res}(\lambda) = - \left(\lambda\!+\!\f\!-\!3\right) Q_4 -
\left(\lambda\!+\!\f\!-\!2\right) P_2(Q_2)
\end{equation}
we find that, for general $n$, the leading coefficient of $-6
Q_8^{res}(\lambda)$ is given by the {\em same} sum as in
\eqref{result-1}. Thus, for the proof of Theorem \ref{Q8-univ} it
suffices to verify

\begin{prop}\label{eval} $-Q_8^{res}(\lambda)^{[4]} = 4! 2^8 w_8$.
\end{prop}

Proposition \ref{eval} is a consequence of the following general
description of the leading coefficients of the families
$P_{2N}(\lambda)$.

\begin{thm}\label{fundamental} For even $n \ge 2$ and $2N \le n$,
the leading coefficient of the degree $N$ polynomial
$$
\left (\lambda\!-\!\f\!+\!1\right) \cdots
\left(\lambda\!-\!\f\!+\!N\right) \T_{2N}(\lambda)
$$
is the multiplication operator by the function
\begin{equation}
\left(v^{-\frac{1}{2}}\right)^{[2N]}.
\end{equation}
\end{thm}

Here we use the notation $f^{[N]}$ for the coefficient of $r^{N}$ in
the Taylor series of $f$. In particular, $v^{[2N]} = v_{2N}$. We
illustrate Theorem \ref{fundamental} by two examples.

\begin{ex} A calculation shows that
\begin{equation}\label{v-4}
\left(v^{-\frac{1}{2}}\right)^{[4]} = -\frac{1}{2} \left( v_4 +
\frac{3}{4} v_2^2 \right).
\end{equation}
Hence
\begin{align*}
P_4(\lambda)^{[2]} & = 2^4 2!
\left(\left(\lambda\!-\!\f\!+\!1\right)\left(\lambda\!-\!\f\!+\!2\right)
\T_4(\lambda)\right)^{[2]} & \mbox{(by definition)} \\
& = 2^4 2! \left(v^{-\frac{1}{2}}\right)^{[4]} & \mbox{(by Theorem \ref{fundamental})} \\
& = -16 v_4 + 12 v_2^2 & \mbox{(by \eqref{v-4})} \\
& = \J^2 + 2|\Rho|^2.
\end{align*}
This result fits with \eqref{P4}.
\end{ex}

\begin{ex} A calculation shows that
\begin{equation}\label{v-6}
\left(v^{-\frac{1}{2}}\right)^{[6]} = -\frac{1}{2}
\left(v_6-\frac{3}{2} v_2 v_4 + \frac{5}{8} v_2^3 \right).
\end{equation}
Hence
\begin{align*}
P_6(\lambda)^{[3]} & = - 2^6 3!
\left(\left(\lambda\!-\!\f\!+\!1\right) \cdots
\left(\lambda\!-\!\f\!+\!3\right)
\T_6(\lambda)\right)^{[3]} & \mbox{(by definition)} \\
& = - 2^6 3! \left(v^{-\frac{1}{2}}\right)^{[6]} & \mbox{(by Theorem \ref{fundamental})} \\
& = 192 \left(v_6-\frac{3}{2} v_2 v_4 + \frac{5}{8} v_2^3 \right) &
\mbox{(by \eqref{v-6}).}
\end{align*}
This result fits with \eqref{P6-leading}.
\end{ex}

Now, for general $N$, combining \eqref{Q-pol} with Theorem
\ref{fundamental} yields
\begin{equation}\label{interm}
Q_{2N}^{res}(\lambda)^{[N]} = -2^{2N} N! \left[
(v^{-\frac{1}{2}})^{[2N]} v_0 + \dots +  (v^{-\frac{1}{2}})^{[0]}
v_{2N} \right].
\end{equation}
But the right-hand side of \eqref{interm} coincides with
$$
-2^{2N} N! \left( v^{-\frac{1}{2}} v \right)^{[2N]} = -2^{2N} N!
(v^{\frac{1}{2}})^{[2N]} = -2^{2N} N! w_{2N}.
$$
Hence we have proved

\begin{prop}\label{leading} For even $n \ge 2$ and $2N \le n$,
$$
Q_{2N}^{res}(\lambda)^{[N]} = -2^{2N} N! w_{2N}.
$$
\end{prop}

In particular, this result proves Conjecture 1.6.2 in \cite{BJ}.

For $N=4$, we obtain Proposition \ref{eval}.

We continue with the

\begin{proof}[Proof of Theorem \ref{fundamental}] A straightforward calculation
shows that the Laplace operator of the metric $g_+ = r^{-2}(dr^2 +
g_r)$ is given by the formula
$$
\Delta_{g_+} = r^2 \frac{\partial^2}{\partial r^2} - (n\!-\!1)r
\frac{\partial }{\partial r} + \frac{1}{2} r^2
\frac{\partial}{\partial r} (\log \det g_r) \frac{\partial}{\partial
r} + r^2 \Delta_{g_r}.
$$
We write the Taylor series of the even function $D(r) = \log \det
g_r$ in the form
$$
D(r) = D^{(0)} + r^2 D^{(1)} + \frac{r^4}{2!} D^{(2)} + \cdots +
\frac{r^n}{(\f)!} D^{(\f)} + \cdots
$$
and expand
$$
\Delta_{g_r} = \Delta^{(0)} + r^2 \Delta^{(1)} + \frac{r^4}{2!}
\Delta^{(2)} + \cdots + \frac{r^n}{(\f)!} \Delta^{(\f)} + \cdots.
$$
The ansatz
$$
u \sim \sum_{N \ge 0} r^{\lambda+2N} \T_{2N}(\lambda)(f), \;
\T_0(\lambda)(f) = f
$$
for solutions of the eigen-equation
$$
-\Delta_{g_+} u = \lambda(n-\lambda) u
$$
leads to the recursive relations
\begin{multline}\label{basic-system}
(\Delta^{(0)} \!+\! (2N\!-\!2\!+\!\lambda) D^{(1)})
\T_{2N-2}(\lambda)(f) +
\dots + \frac{1}{(N\!-\!1)!} (\Delta^{(N-1)} \!+\! \lambda D^{(N)}) \T_0(\lambda)(f) \\
= -2N(2\lambda\!-\!n\!+\!2N) \T_{2N}(\lambda)(f).
\end{multline}
For $N=3$, the latter formula yields \eqref{p6-form}. Let
$\omega_{2N}$ be the leading coefficient of the polynomial
$$
\left(\lambda\!-\!\f\!+\!N\right) \cdots
\left(\lambda\!-\!\f\!+\!1\right) \T_{2N}(\lambda)(f).
$$
We also set $\omega_0=f$. Eq.~\eqref{basic-system} shows that the
coefficients $\omega_{2N}$ are recursively determined by the system
of relations
\begin{equation}\label{lead-system}
D^{(1)} \omega_{2N-2} + \cdots + \frac{1}{(N\!-\!1)!} D^{(N)}
\omega_0 = -4N \omega_{2N}, \; N \ge 1.
\end{equation}
Now
$$
-\sum_{N \ge 1} 4N \omega_{2N} r^{2N-2} = -2\frac{1}{r}
\frac{\partial}{\partial r} \left( \sum_{N \ge 0} \omega_{2N} r^{2N}
\right)
$$
and
\begin{multline*}
\sum_{N \ge 1} \left(D^{(1)} \omega_{2N-2} + \cdots +
\frac{1}{(N\!-\!1)!} D^{(N)} \omega_0 \right) r^{2N-2} \\
= \frac{1}{2} \frac{1}{r} \frac{\partial}{\partial r} \left(D^{(0)}
+ r^2 D^{(1)} + \cdots + \frac{r^n}{(\f)!} D^{(\f)} + \cdots \right)
\left( \sum_{N \ge 0} \omega_{2N} r^{2N} \right) \\
= \frac{1}{r} \frac{\partial}{\partial r} \log v(r) \left( \sum_{N
\ge 0} \omega_{2N} r^{2N} \right).
\end{multline*}
Here we suppressed that, for general metrics $g$, the quantities are
only well-defined for $2N \le n$. If necessary, the identities are
to be interpreted as those for finite sums. Now the assertion
follows from the fact that the function $\psi = v^{-\frac{1}{2}}$
satisfies the differential equation
$$
-2 \frac{\partial}{\partial r} \psi = \frac{\partial}{\partial r}
(\log v) \psi.
$$
The proof is complete. \end{proof}

Theorem \ref{fundamental} also implies the following result. Let $n$
be even and $2N \le n$.

\begin{prop}\label{V-pol-gen} On any manifold of even dimension $n$,
the leading coefficient of the degree $N$ polynomial
\begin{equation}\label{V-pol-def}
\V_{2N}(\lambda) = \left[\left(\lambda\!-\!\f\!+\!1\right) \cdots
\left(\lambda\!-\!\f\!+\!N\right) \right] \sum_{j=0}^N (2N\!+\!2j)
\T^*_{2j}(\lambda)(v_{2N-2j})
\end{equation}
vanishes, i.e.,
\begin{equation}
\V_{2N}(\lambda)^{[N]} = 0.
\end{equation}
\end{prop}

\begin{proof} As above, the following identities are to be interpreted
as relations for terminating sums (if necessary). By Theorem
\ref{fundamental}, the assertion is equivalent to
\begin{equation}\label{van-gen}
\sum_{N \ge 0} \left( \sum_{j=0}^N (2N\!+\!2j)
\left(v^{-\frac{1}{2}}\right)^{[2j]} v^{[2N-2j]} \right) r^{2N} = 0.
\end{equation}
Now we have
\begin{equation*}
\sum_{N \ge 0} 2N \left( \sum_{j=0}^N
\left(v^{-\frac{1}{2}}\right)^{[2j]} v^{[2N-2j]} \right) r^{2N}  = r
\frac{\partial}{\partial r} (v^{-\frac{1}{2}} v) = \frac{1}{2} r
v^{-\frac{1}{2}} \frac{\partial v}{\partial r}
\end{equation*}
and
\begin{equation*}
\sum_{N \ge 0} \left( \sum_{j=0}^N 2j
\left(v^{-\frac{1}{2}}\right)^{[2j]} v^{[2N-2j]} \right) r^{2N} \\
= r \frac{\partial}{\partial r} (v^{-\frac{1}{2}}) v = -\frac{1}{2}
r v^{-\frac{1}{2}} \frac{\partial v}{\partial r}.
\end{equation*}
Summing both relations proves \eqref{van-gen}.
\end{proof}

Similarly as in Section \ref{proof1}, one can use the vanishing result
$$
\V_8(\lambda)^{[4]} = 0
$$
to express $Q_8^{res}(\lambda)^{[4]}$ in terms of $P_2(\lambda)$,
$P_4(\lambda)$, $P_6(\lambda)$ and $v_2,\dots,v_8$. This gives another
proof of Theorem \ref{Q8-univ}.

By Proposition \ref{V-pol-gen}, the degree of the polynomial
$\V_{2N}(\lambda)$ is $N\!-\!1$. In contrast to the critical case,
the polynomial $\V_{2N}(\lambda)$ does not vanish identically if $2N
< n$. In the following example we make $\V_2$ and $\V_4$ explicit.
For a further discussion of $\V_{2N}(\lambda)$ see Section
\ref{final-co}.

\begin{ex}\label{V-pol} We have $\V_2(\lambda) = \left(\f\!-\!1\right) Q_2$,
and the linear polynomial $\V_4(\lambda)$ equals
$$
\V_4(\lambda) = \frac{1}{4} \left(\f\!-\!2\right) \left[ -
\left(\lambda\!-\!\f\!+\!2\right) (Q_4+P_2(Q_2)) + Q_4 \right].
$$
These formulas show that
$$
\V_2(\lambda) = \left(\f\!-\!1\right) \QR_2^{res}(\lambda\!-\!n\!+\!2)
\quad \mbox{and} \quad \V_4(\lambda) = \frac{1}{4}\left(\f\!-\!2\right)
\QR_4^{res}(\lambda\!-\!n\!+\!4)
$$
using $\QR_2^{res}(\lambda) = Q_2$ and \eqref{Q4-gen}.
\end{ex}

\section{Recursive formulas for $Q_4$, $Q_6$ and $P_4$, $P_6$}
\label{final}

Theorems \ref{main} -- \ref{reduced} are analogs of similar results
for the $Q$-curvatures $Q_4$ and $Q_6$. Here we sketch proofs of
these formulas. In addition, we display universal recursive formulas
for $P_4$ and $P_6$. In combination with Theorems \ref{main} --
\ref{reduced}, these results can be used to derive more explicit
formulas for $Q_8$.

We start with the discussion of $Q_4$ and $Q_6$.

\begin{prop}\label{Q4-both} On manifolds of dimension $n \ge 4$,
\begin{equation}\label{rec-Q4a}
Q_4 = -P_2(Q_2) - Q_2^2 + 2! 2^3 v_4.
\end{equation}
This formula is equivalent to
\begin{equation}\label{rec-Q4b}
Q_4 = -P_2(Q_2) + 2! 2^4 w_4
\end{equation}
with
$$
8 w_4 = 4v_4 - v_2^2.
$$
In dimension $n=4$, the reduced form of \eqref{rec-Q4a} reads
\begin{equation}\label{rec-Q4c}
Q_4 = -P_2^0(Q_2) + 2! 2^3 v_4.
\end{equation}
In particular,
$$
\int_{M^4} Q_4 vol = 2! 2^3 \int_{M^4} v_4 vol.
$$
\end{prop}

We recall that $8 v_4 = \J^2 - |\Rho|^2$ (see \eqref{v24}). The
assertions are simple consequences of \eqref{yamabe} and
\eqref{q4-gen}.

\begin{prop}\label{Q6-both} On manifolds of dimension $n \ge 6$,
\begin{equation}\label{rec-Q6a}
Q_6 = \left[-2P_2(Q_4) + 2 P_4(Q_2) - 3 P_2^2(Q_2)\right] - 6 \left[
Q_4 + P_2(Q_2) \right] Q_2  - 2!3!2^5 v_6.
\end{equation}
This formula is equivalent to
\begin{equation}\label{rec-Q6b}
Q_6 = \left[-2P_2(Q_4) + 2 P_4(Q_2) - 3 P_2^2(Q_2)\right] - 2!3!2^6
w_6
\end{equation}
with
$$
16 w_6 = 8v_6 - 4v_4v_2 + v_2^3.
$$
In dimension $n=6$, the reduced form of \eqref{rec-Q6a} reads
\begin{equation}\label{rec-Q6c}
Q_6 = \left[-2P_2^0(Q_4) + 2 P_4^0(Q_2) - 3 P_2^0 P_2 (Q_2)\right] -
2!3!2^5 v_6.
\end{equation}
In particular,
$$
\int_{M^6} Q_6 vol = 2! 3! 2^5 \int_{M^6} v_6 vol.
$$
\end{prop}

We recall that
\begin{equation}\label{v6}
8 v_6 = - \tr (\wedge^3 \Rho) - \frac{1}{3(n\!-\!4)} (\B,\Rho),
\end{equation}
where
\begin{equation}\label{bach}
\B_{ij} = \Delta (\Rho)_{ij} - \nabla^k \nabla_j (\Rho)_{ik} +
\Rho^{kl} \W_{kijl}
\end{equation}
generalizes the Bach tensor. For a proof of \eqref{v6} see Theorem
6.9.2 in \cite{juhl-book}.

\begin{proof} Let $n=6$. We first sketch a proof of \eqref{rec-Q6a}
along the same lines of arguments as in Section \ref{proof1}. By
Theorem 6.11.9 in \cite{juhl-book}, we have
$$
Q_6^{res}(\lambda) = \frac{1}{2} \lambda
(\lambda\!-\!1)(\lambda\!-\!2) Q_6 + \lambda^2(\lambda\!-\!1) P_2
\left(Q_4+ \frac{3}{2} P_2(Q_2)\right) - \lambda^2(\lambda\!-\!2)
P_4(Q_2).
$$
This formula is a consequence of the $3$ factorizations
$$
Q_6^{res}(2) = - P_2(Q_4^{res}(2)), \; Q_6^{res}(1) = P_4
(Q_2^{res}(1)), \; Q_6^{res}(0) = P_6 (1) = 0
$$
and $\dot{Q}_6^{res}(0) = Q_6$. It follows that the coefficient of
$\lambda^3$ in $2Q_6^{res}(\lambda)$ is given by
\begin{equation}\label{top-1}
Q_6 + 2P_2(Q_4) - 2P_4(Q_2) + 3P_2^2(Q_2).
\end{equation}
On the other hand, by definition, $Q_6^{res}(\lambda)$ equals
\begin{equation}\label{q6-rep1}
P_6^*(\lambda)(v_0) - 12 \lambda P_4^*(\lambda)(v_2) + 2^4 3!
\lambda(\lambda\!-\!1) P_2^*(\lambda)(v_4) - 2^6 3!
\lambda(\lambda\!-\!1)(\lambda\!-\!2) v_6.
\end{equation}
In this sum, the term $P_6^*(\lambda)(1)$ can be written as a linear
combination of the other three terms. We express this fact as the
vanishing result
\begin{equation}\label{V6-van}
\V_6(\lambda) = 0.
\end{equation}
In fact, \eqref{V6-van} follows from
$$
\V_6(0) = \V_6(1) = \V_6(2) = 0 \quad \mbox{and} \quad \dot{\V}_6(0)
= 0.
$$
In Section \ref{proof1}, we have seen that the first and the last
properties are special cases of a general result, and we proved
$\V_6(1) = \V_6(2) = 0$ (for $n=6$).\footnote{For an alternative
proof of the identity \eqref{q6-rep2} we refer to Lemma 6.11.10 in
\cite{juhl-book}.} Now \eqref{V6-van} is equivalent to the identity
\begin{equation}\label{q6-rep2}
P_6^*(\lambda)(1) = 10 \lambda P_4^*(\lambda)(v_2) - 2^6
\lambda(\lambda\!-\!1) P_2^*(\lambda)(v_4) + 2^6 3
\lambda(\lambda\!-\!1)(\lambda\!-\!2) v_6.
\end{equation}
Thus, \eqref{q6-rep1} implies
\begin{equation}
Q_6^{res}(\lambda) = -2\lambda P_4^*(\lambda)(v_2) + 2^5
\lambda(\lambda\!-\!1) P_2^*(\lambda)(v_4) - 2^6 3
\lambda(\lambda\!-\!1)(\lambda\!-\!2) v_6,
\end{equation}
and using \eqref{P4} we find
\begin{equation}\label{top-2}
2 Q_6^{res}(\lambda)^{[3]} = - 6 \J^3 + 12 \J |\Rho|^2 - 2^6 3! v_6.
\end{equation}
Now a comparison of \eqref{top-1} and \eqref{top-2} implies
$$
Q_6 + 2P_2(Q_4) - 2P_4(Q_2) + 3P_2^2(Q_2) = - 6 (\J^2 - 2|\Rho|^2)
\J - 2^6 3! v_6.
$$
We complete the proof of \eqref{rec-Q6a} by rewriting this identity
using \eqref{rec-Q4a} (in dimension $n=6$).

The identities \eqref{rec-Q6b} and \eqref{rec-Q6c} follow from
\eqref{rec-Q6a} by straightforward calculations.

In \cite{juhl-book}, we gave a proof that \eqref{rec-Q6a} remains
valid in all dimensions $n \ge 6$. It rests on an explicit formula
for $Q_6$ which follows from a combination the relation
$$
P_6\left(\f\!-\!3\right)(1) = -\left(\f\!-\!3\right) Q_6
$$
with the formula \eqref{p6-form}. Next, we present an alternative
proof of the universality of \eqref{rec-Q6b} along the lines of
Section \ref{universal}. For $n \ge 6$, \eqref{Q6-gen} and
\eqref{Q4-gen} show that the leading coefficient of
$2Q_6^{res}(\lambda)$ is still given by the linear combination
\eqref{top-1}. On the other hand, Proposition \ref{leading} implies
$$
2Q_6^{res}(\lambda)^{[3]} = -2!3!2^6 w_6.
$$
The comparison of both results proves the universality of
\eqref{rec-Q6b}. Now \eqref{rec-Q6a} follows by direct calculation
using Proposition \ref{Q4-both}. Note that the latter arguments do
not require to know explicit expressions for the quantities
involved.
\end{proof}

We stress that the multiplicities in \eqref{rec-Q4a} and
\eqref{rec-Q6a} again are given by the general rule \eqref{mult}.

The following results describe the non-constant parts of the
GJMS-operators $P_4$ and $P_6$. Of course, their constant terms are
given by the corresponding $Q$-curvatures.

\begin{prop}\label{P4-univ} On manifolds of dimension $n \ge 4$,
\begin{equation}\label{rec-P4}
P_4^0 = (P_2^2)^0 - 4 \delta (\Rho d).
\end{equation}
\end{prop}

The result follows by a calculation from \eqref{pan-op}.

\begin{prop}\label{P6-univ} On manifolds of dimension $n \ge 6$,
\begin{equation}\label{rec-P6}
P_6^0 = \left[ 2 P_2 P_4 + 2 P_4 P_2 - 3 P_2^3 \right]^0 - 48
\delta(\Rho^2 d) - \frac{16}{n\!-\!4} \delta(\B d).
\end{equation}
\end{prop}

Proposition \ref{P6-univ} follows from \eqref{rec-Q6a} by
infinitesimal conformal variation. For the details we refer to
Section 6.12 of \cite{juhl-book}.

\section{Final comments}\label{final-co}

Among all $Q$-curvatures of a manifold of even dimension $n$, the
critical $Q$-curvature $Q_n$ is distinguished by the property that its
behaviour under conformal changes of the metric is governed by the
{\em linear} differential operator $P_n$. More precisely, the pair
$(P_n,Q_n)$ satisfies the fundamental identity
$$
e^{n\varphi} Q_n(e^{2\varphi}g) = Q_n(g) + (-1)^\f P_n(g)(\varphi)
$$
for all $\varphi \in C^\infty(M)$. Proposition \ref{P6-univ} shows
that, up to a second-order operator, the critical $P_6$ can be
written as a linear combination of compositions of lower order
GJMS-operators. Moreover, the multiplicities of the compositions in
that sum are related to the multiplicities of corresponding terms in
the recursive formula \eqref{rec-Q6a}. Along the same line, Theorem
\ref{main} suggests to expect an analogous recursive formula for the
critical GJMS-operator $P_8$. Theorem 11.1 in \cite{juhl-power}
establishes the conformal covariance of such a natural candidate for
the GJMS-operator $P_8$. It remains an open problem to prove that
this operator actually coincides with $P_8$.

It is well-known (\cite{bran-2}, Corollary 1.5) that the contribution
to $Q_{2N}$ which involves the maximal number of derivatives is given
by $(-1)^{n-1} \Delta^{N-1} \J$. This is obvious for $Q_4$ and can be
reproduced for $Q_6$ and $Q_8$ by using \eqref{rec-Q6a} and
\eqref{Q8-main}, respectively. Indeed, the latter facts are special
cases of the summation formula
$$
\sum_{|I|=N} m_I = 0
$$
(see Lemma 2.1 in \cite{juhl-power}).

A full comparison of \eqref{Q8-main} (in general dimensions) with the
formula of Gover and Peterson for $Q_8$ (in general dimensions) (see
Figure 5 in \cite{G-P}) remains a challenge. As an example, we
consider the contribution of $(\Delta \Rho,\Delta \Rho)$ to $Q_8$.  By
\cite{G-P}, this term has the coefficient
\begin{equation}\label{coeff} -12
\frac{n^2\!-\!4n\!+\!8}{(n\!-\!4)^2} = -12 \left( 1 +
\frac{4}{n\!-\!4} + \frac{8}{(n\!-\!4)^2} \right).
\end{equation}
The result is confirmed by \eqref{Q8-main}. In fact, the term $(\Delta
\Rho,\Delta \Rho)$ contributes to $Q_8$ only through
$$
-12 P_2^2(Q_4), \quad -9 P_4(Q_4), \quad -3P_2(Q_6) \quad \mbox{and}
\quad 3!4!2^7 v_8.
$$
The first two terms yield
$$
48 (\Delta \Rho,\Delta \Rho) \quad \mbox{and} \quad -36 (\Delta
\Rho,\Delta \Rho).
$$
\eqref{rec-Q6a} and \eqref{v6} show that the third terms contributes
by
$$
-24 (\Delta \Rho,\Delta \Rho) - \frac{48}{n\!-\!4} (\Delta
\Rho,\Delta \Rho).
$$
Finally, Graham's formula for $v_8$ and \eqref{bach} show that the
last term contributes by\footnote{Here we correct a misprint in
equation (2.23) of \cite{G-ext}: the term $\tr (\Omega^{(1)})^2$ is
to be replaced by $\tr ((\Omega^{(1)})^2)$. The tensor $\Omega^{(1)}
= \frac{\B}{4-n}$ is the first extended obstruction tensor.}
$$
-\frac{96}{(n\!-\!4)^2} (\Delta \Rho,\Delta \Rho).
$$
Summarizing, we find \eqref{coeff}.

A version of \eqref{Q8-gen} holds true for all $Q$-curvature
polynomials $Q_{2N}^{res}(\lambda)$ with $2 \le 2N \le n$. It
follows that the leading coefficient of $Q_{2N}^{res}(\lambda)$ can
be written as the product of $(-1)^{N-1}(N\!-\!1)!$ and a linear
combination of the form
$$
\sum_{|I|+a=N} \mu_{(I,a)} P_{2I}(Q_{2a}) = Q_{2N} + \mbox{terms
involving lower order $Q$-curvatures}
$$
with certain coefficients $\mu_I \in \r$, $|I|=N$. On the other
hand, Proposition \ref{leading} implies that this sum coincides with
$-2^{2N} N! w_{2N}$. Hence we obtain an identity of the form
$$
\sum_{|I|+a=N} \mu_{(I,a)} P_{2I}(Q_{2a}) = (-1)^N N!(N\!-\!1)!
2^{2N} w_{2N}.
$$
This proves Conjecture 9.2 in \cite{juhl-power}, up to the algebraic
problem to establish the identifications
\begin{equation}\label{id-coeff}
(-1)^a m_{(I,a)} = \mu_{(I,a)} \quad \mbox{for all $(I,a)$}.
\end{equation}
For small $|I|+a$, the relations \eqref{id-coeff} follow either by
the evaluation of the algorithm which generates the formulas for
$Q$-curvature polynomials in terms of $Q$-curvatures and
GJMS-operators, or from the validity of Conjecture 9.2 on round
spheres and pseudo-spheres (proved in \cite{JK}).

Although these arguments suffice to prove Conjecture 9.2 also for,
say, $Q_{10}$, in the present paper we have restricted the attention
to $Q_8$ since only in that case we gain a complete understanding of
$Q_8$ in terms of the metric. In fact, a fully explicit formula in
terms of $\Rho$ and Graham's first two extended obstruction tensors
$\Omega^{(1)}$, $\Omega^{(2)}$ (see \cite{G-ext}) follows by
combining \eqref{Q8-main} with the formulas displayed in Section
\ref{final}. Future applications will shown to which extent such
explicit versions are of interest.

The polynomial $\V_{2N}(\lambda)$ (see \eqref{V-pol-def}) seems to
be related to the $Q$-curvature polynomial $Q_{2N}^{res}(\lambda)$
by the formula
\begin{equation}\label{VQ}
2^{2N-2} (N\!-\!1)! \V_{2N}(\lambda) = \left(\f\!-\!N\right)
\QR_{2N}^{res}(\lambda\!-\!n\!+\!2N).
\end{equation}
In more explicit terms, \eqref{VQ} states the equality
$$
(\lambda\!-\!n\!+\!2N) \sum_{j=0}^N (2N\!+\!2j)
\T_{2j}^*(\lambda)(v_{2N-2j}) = -2N(n\!-\!2N) \sum_{j=0}^N
\T_{2j}^*(\lambda)(v_{2N-2j})
$$
of rational functions in $\lambda$. The special cases
$N=1$ and $N=2$ of \eqref{VQ} appear in Example \ref{V-pol}. For
$N=3$ and $N=4$, the relation \eqref{VQ} can be proved by direct
calculations, too. In fact, by Proposition \ref{V-pol-gen} the
polynomial $\V_{2N}(\lambda)$ has degree $N\!-\!1$, and thus it
suffices to verify that $\V_{2N}(\lambda)$ satisfies the $N$
factorization identities which correspond to those of
$Q_{2N}^{res}(\lambda)$ by Proposition \ref{q-factor}. In
particular, for $\V_6(\lambda)$, these state that
\begin{align*}
32 \V_6 \left(\f-3\right) & = - P_6 (1), \\
32 (n-2) \V_6 \left(\f-2\right) & = (n-6) P_4 \V_2 \left(\f+2\right), \\
32 (n-4) \V_6 \left(\f-1\right) & = (n-6) P_2 \V_4
\left(\f+1\right).
\end{align*}
In the critical case $2N=n$, \eqref{VQ} would imply that
$\V_n(\lambda)=0$. This is the assertion of Conjecture 6.11.2 in
\cite{juhl-book}. Since $\QR_{2N}^{res}(\lambda)$ has degree
$N\!-\!1$, \eqref{VQ} would also imply that $\V_{2N}(\lambda)^{[N]}
= 0$, i.e., Proposition \ref{V-pol-gen}.

Finally, we note that alternative universal recursive formulas for
$Q$-curvatures can be derived by using some of the additional
identities which are satisfied by the $Q$-curvature polynomials (see
\cite{FJ}).


\end{document}